\def\ara{\mathop{\rm ara}\nolimits}
\def\Ass{\mathop{\rm Ass}\nolimits}
\def\Assh{\mathop{\rm Assh}\nolimits}
\def\codim{\mathop{\rm codim}\nolimits}

\def\height{\mathop{\rm height}\nolimits}
\def\Hom{\mathop{\rm Hom}\nolimits}

\def\im{\mathop{\rm im}\nolimits}
\def\ker{\mathop{\rm ker}\nolimits}

\def\Ext{\mathop{\rm Ext}\nolimits}
\def\InjH{\mathop{\rm E}\nolimits}
\def\LCMo{\mathop{\rm H}\nolimits}

\def\Spec{\mathop{\rm Spec}\nolimits}
\def\Supp{\mathop{\rm Supp}\nolimits}

\def\Naturalsign{{\rm l\kern-.23em N}}
\font \normal=cmr10 scaled \magstep0 \font \mittel=cmr10 scaled \magstep1 \font \gross=cmr10 scaled \magstep5
\input amssym.def
\input amssym.tex
{
\parindent=0pt
\gross
Matlis duals of top Local Cohomology Modules
\normal
\bigskip
\bigskip
Michael Hellus, michael.hellus@math.uni-leipzig.de
\par
J\"urgen St\"uckrad, juergen.stueckrad@math.uni-leipzig.de
\footnote{$ $}{2000 MSC: 13D45, 13C05.} \footnote{$ $}{Key words:
Local cohomology, Matlis duality, associated prime ideals.}
\bigskip
Abstract: In the first section of this paper we present
generalizations of known results on the set of associated primes of
Matlis duals of local cohomology modules; we prove these
generalizations by using a new technique. In section 2 we compute
the set of associated primes of the Matlis dual of $\LCMo
^{d-1}_J(R)$, where $R$ is a $d$-dimensional local ring and
$J\subseteq R$ an ideal such that $\dim (R/J)=1$ and $\LCMo
^d_J(R)=0$.
\bigskip
\mittel
0. Introduction and Motivation
\normal
\bigskip
}
In algebraic geometry one investigates the number of equations needed to cut out a
given variety; this number is called the arithmetic rank $\ara$; thus for a given ideal $I$ in a
noetherian ring $R$ we have
$$\ara (I)=\min \{ h\in \Naturalsign \vert \exists
r_1,\dots ,r_h\in I: \sqrt I=\sqrt {(r_1,\dots ,r_h)R}\} .$$ One
always has $\ara (I)\geq \codim (I)$. We say $I$ (or the variety
given by $I$) is a (set-theoretic) complete intersection if $\ara
(I)=\codim (I)$ holds. For an $R$-module $M$, we denote by $\LCMo
^i_I(M)$ the i-th local cohomology module of $M$ supported in $I$.
Let $n$ be a natural number; then $\ara (I)\leq n$ implies $0=\LCMo
^{n+1}_I(R)=\LCMo ^{n+2}_I(R)=\dots $, but the converse implication
does not hold in general (see [9, example 1.1] for a
counterexample). If $I$ is a complete intersection in a
Cohen-Macaulay ring $R$ with $\ara (I)=:h$ then $\LCMo ^l_I(R)\neq
0\iff l=h$ holds. More precisely, at least if $R$ is local, the case
$I$ is a complete intersection can be characterized as follows:
Assume that $\LCMo ^l_I(R)\neq 0\iff l=h$ holds. Let $r_1,\dots
,r_h\in I$ be an $R$-regular sequence; then
$$\sqrt
I=\sqrt {(r_1,\dots ,r_h)R}\iff r_1,\dots ,r_h\hbox { is a regular
sequence on }D(\LCMo ^h_I(R))$$ holds (for a proof see [11, section
0] and for related remarks see [9, section 3]). Here $D(\_ ):=D_R(\_
)$ stands for the Matlis dual functor with respect to the local ring
$R$ (i. e. $D(\_ )=\Hom _R(\_ \ ,\InjH _R(R/\goth m))$, where $\InjH
_R(M)$ stands for an $R$-injective hull of an $R$-module $M$) and
the condition $r_1,\dots ,r_h$ is a regular sequence on $D(\LCMo
^h_I(R)$ means that $D(\LCMo ^h_I(R))/(r_1,\dots ,r_h)D(\LCMo
^h_I(R))\neq 0$ holds and that every $r_i$ operates injectively on
$D(\LCMo ^h_I(R))/(r_1,\dots ,r_{i-1})D(\LCMo ^h_I(R))=D(\LCMo
^{h-i+1}_{(r_i,\dots ,r_h)R}(R/(r_1,\dots ,r_{i-1})R))$. Thus if one
is interested in the arithmetic rank or in complete intersections
one is naturally lead to consider the set of associated primes of
Matlis duals of top local cohomology modules. Similar questions were
examined in [11] and [9]. In [11, section 1] it was conjectured that
for every noetherian local ring $(R,\goth m)$ the equality
$$\Ass _R(D(\LCMo ^i_{(x_1,\dots ,x_i)R}(R)))=\{ \goth p\in \Spec (R)\vert
\LCMo ^i_{(x_1,\dots ,x_i)R}(R/\goth p)\neq 0\} \eqno {(*)}$$ holds
for any $i\geq 1$, $x_1,\dots ,x_i\in R$; there it was shown ([11,
theorem 1.1]) that this conjecture (*) is equivalent to: If
$(R,\goth m)$ is a noetherian local ring, $i\geq 1$ and $x_1,\dots
,x_i\in R$, the set $Y:=\Ass _R(D(\LCMo ^i_{(x_1,\dots ,x_i)R}(R)))$
is stable under generalization, i. e.
$$\goth p,\goth q\in \Spec (R), \goth p\subseteq \goth q,\goth
q\in Y\Rightarrow \goth p\in Y.$$
holds. It seems reasonable to conjecture more, namely: In the situation above all primes
$\goth p$ maximal in $Y$ have the same dimension, namely $\dim (R/\goth
p)=i$. We will refer to this conjecture by (+).
\par
In [11, theorem 2.2.1] it was shown that for every noetherian local
ring $R$ containing a field
$$\{ \goth p\in \Spec (R)\vert x_1,\dots ,x_i\hbox { is part of a system of
parameters of }R/\goth p\} \subseteq \Ass _R(D(\LCMo ^i_{(x_1,\dots
,x_i)R}(R)))$$ holds for any $i\geq 1$, $x_1,\dots ,x_i\in R$;
secondly it was shown (see [11], corollary 2.2.2) that for every
such ring $R$ and every $x\in R$ one has
$$\Ass _R(D(\LCMo ^1_{xR}(R)))=\Spec (R)\setminus {\frak V}(x).$$
And thirdly it is known ([9, theorem 1.1]) that for an ideal $I$ in
a noetherian local ring $R$ such that $0=\LCMo ^2_I(R)=\LCMo
^3_I(R)=\dots $ holds, one has
$$\ara (I)\leq 1\iff \Ass _R(D(\LCMo ^1_I(R)))\hbox { satisfies prime avoidance}$$
(see [9, theorem 1.1 (i) and (ii)] for details).
\par
In section 1 of this work we will generalize the first two of the
last three mentioned results to the class of all noetherian local
rings; this will be done by methods different to the ones used in
[11]; in this section the crucial point is lemma 1.1. In addition we
will prove some more facts related to these results, e. g. a version
of the first of the above mentioned results for local cohomology of
an $R$-module $M$ (instead of $R$).
\par
We would like to point out that Brodmann and Huneke proved lemma 1.6
(that will be needed in section 2) and used it to give a short proof
of the Hartshorne-Lichtenbaum vanishing theorem (see [3]).
\par
In section 2 we consider a $d$-dimensional local ring, an ideal $J$
of $R$ such that $\dim (R/J)=1$ and $\LCMo ^d_J(R)=0$. We are
interested in $\Ass _R(D(\LCMo ^{d-1}_J(R)))$. We obtain partial
results in the general case; in the case where $R$ is complete we
are able to completely compute this set:
$$\Ass _R(D(\LCMo ^{d-1}_J(R))=\{ \goth p\in \Spec (R)\vert \dim (R/\goth
p)=d-1, \dim (R/(\goth p+J)=0\} \cup \Assh R.$$
Thus $\Ass _R(D(\LCMo ^{d-1}_J(R))$ is stable under generalization. $\Assh (M)$ denotes the set of associated primes of $M$ of highest dimension
($M$ any $R$-module).
\parindent=0pt
\bigskip
\mittel 1. Preliminaries
\normal
\bigskip
{\bf 1.1 Lemma}
\par
Let $R$ be a ring, $x,y\in R$ and $U$ an $R$-submodule of $R_x$ such that $\im
\iota _x\subseteq U$, where $\iota _x:R\to R_x$ is the canonical map. Let $S:=
\im \iota _y\subseteq R_y$. There exists an $R$-epimorphism
$$R_x/U\to R_{xy}/(S_x+U_y).$$
Proof: Let $V:=S_x+U_y\subseteq R_{xy}$ and let $(b_1,b_2,\dots )\in
R^{\Naturalsign ^+}$ be an infinite sequence. For $i\in \Naturalsign $ we set
$$\rho _i:=\sum _{j=1}^i{b_j\over x^{i-j+1}y^j}+V\in R_{xy}/V\ \ (i\in
\Naturalsign ).$$
We calculate
$$\eqalign {x\rho _{i+1}-\rho _i&=(\sum _{j=1}^{i+1}{xb_j\over
x^{i-j+2}y^j}+V)-(\sum _{j=1}^i{b_j\over x^{i-j+1}y^j}+V)\cr &={b_{i+1}\over
y^{i+1}}+V\cr &=0,\cr }$$
because ${b_{i+1}\over y^{i+1}}\in (\im \iota _x)_y\subseteq U_y\subseteq
V$. Thus we have $x\rho _{i+1}=\rho _i$ for all $i\in \Naturalsign $ and so we
get a map $\varphi :R_x\to R_{xy}/V$ given by
$${r\over x^i}\mapsto r\rho _i\ \ (r\in R,i\in \Naturalsign ).$$
It is easy to see that $\varphi $ is $R$-linear. Let $u\in U$ be
arbitrary. There are $r\in R$ and $i\in \Naturalsign $ such that $u={r\over
x^i}$. We have
$$\varphi (u)=r\rho _i=\sum _{j=1}^i{rb_j\over x^{i-j+1}}+V=u\sum
_{j=1}^i{x^{j-1}b_j\over y^j}+V=0,$$
because $u\sum _{j=1}^i{x^{j-1}b_j\over y^j}\in U_y\subseteq V$. This implies
$U\subseteq \ker (\varphi )$ and hence we get an induced $R$-homomorphism
$f:R_x/U\to R_{xy}/V$. The set $\{ {1\over x^i}+U\vert i\in \Naturalsign ^+\}
$ is a generating set for $R_x/U$ and so we have
$$f\hbox { is surjective }\iff \varphi \hbox{ is surjective }\iff \{ \rho _1,\rho
_2, \dots \} \hbox { generates } R_{xy}/V.$$
The set $\{ {1\over x^iy^j}+V\vert i,j\in \Naturalsign ^+\} $ generates
$R_{xy}/V$. For $i\in \Naturalsign ^+$ we set
$${\cal B}_i:=\pmatrix {&b_1&b_2&b_3&\ldots &b_i\cr &b_2&b_3&b_4&\ldots
&b_{i+1}\cr &\vdots &\vdots &\vdots & &\vdots \cr &b_i&b_{i+1}&b_{i+2}&\ldots
&b_{2i-1}\cr }$$
Then we have for $i\in \Naturalsign ^+$:
$$(\rho _i,y\rho _{i+1},\dots ,y^{i-1}\rho _{2i-1})^T={\cal B}_i({1\over
x^iy}+V,{1\over x^{i-1}y^2}+V,\dots ,{1\over xy^i}+V)^T.$$
If we choose $b_1,b_2,\dots \in R$ in such a way that $\det {\cal B}_i\in R^*$ for
all $i\in \Naturalsign ^+$ (which is possible, $\cal B$ consists only of ones
and zeroes), then $\{ \rho _1,\rho _2,
\dots \}$ is generating $R_{xy}/V$.
\bigskip
From now on we assume $R$ is noetherian and we can use Cech
cohomology.
\bigskip
{\bf 1.2 Theorem}
\par
Let $R$ be a noetherian ring, $M$ an $R$-module, $m\in
\Naturalsign ^+, n\in \Naturalsign ,x_1,\dots ,x_m,y_1,\dots ,y_n\in R$. Then
there exists an $R$-epimorphism
$$\LCMo ^m_{(x_1,\dots ,x_m)R}(M)\to \LCMo ^{m+n}_{(x_1,\dots ,x_m,y_1,\dots
,y_n)R}(M).$$
Proof: Obviously it suffices to prove the statement for the case $M=R$. Using Cech cohomology to compute both local cohomology modules the
statement follows immediately from lemma 1.1 by induction on $n$.
\bigskip
{\bf 1.3 Theorem}
\par
Let $(R,\goth m)$ be a noetherian local ring, $m\in
\Naturalsign ^+, x_1,\dots ,x_m\in \goth m$ and $M$ a finitely generated
$R$-module. Then the following statements hold:
\par
(i) $\dim (M/\goth pM)\geq m$ for every $\goth p\in \Ass _R(D(\LCMo
^m_{(x_1,\dots ,x_m)R}(M)))$.
\par
(ii) $\{ \goth p\in \Supp _R(M)\vert x_1,\dots ,x_m\hbox { is part of a system
of parameters of }R/\goth p\} \subseteq \Ass _R(D(\LCMo ^m_{(x_1,\dots
,x_m)R}(M)))$.
\par
(iii) $\Ass _R(D(\LCMo ^1_{xR}(R)))=\Spec (R)\setminus {\frak V}(x)$ for every
$x\in R$.
\par
(iv) If $x_1,\dots ,x_m$ is part of a system of parameters of $M$, we have
$\Assh (M)\subseteq \Ass _R(D(\LCMo ^m_{(x_1,\dots ,x_m)R}(M)))$; furthermore,
in case $m=\dim (M)$ equality holds: $\Ass _R(D(\LCMo ^{\dim (M)}_\goth
m(M)))=\Assh (M)$.
\par
(v) If $R$ is complete, $\goth p\in \Supp _R(M)$ and $\dim (R/\goth p)=m$, the
equivalence
$$\goth p\in \Ass _R(D(\LCMo ^m_{(x_1,\dots ,x_m)R}(M)))\iff x_1,\dots
,x_m\hbox { is a system of parameters of }R/\goth p$$
holds.
\par
Proof:
We set $I:=(x_1,\dots ,x_m)R$.
\par
(i) Let $\goth p\in \Ass _R(D(\LCMo ^m_I(M)))$. We conclude
$$0\neq \Hom _R(R/\goth p,D(\LCMo ^m_I(M)))=D(\LCMo ^m_I(M)\otimes _R(R/\goth
p))=D(\LCMo ^m_I(M/\goth pM)).$$
Thus we have $\LCMo ^m_I(M/\goth pM)\neq 0$ and statement (i) follows.
\par
(ii) Let $\goth p\in \Supp _R(M)$ such that $x_1,\dots ,x_m$ is part of a
system of parameters of $R/\goth p$. By completing $x_1,\dots ,x_m$ to a
system of parameters of $M/\goth pM$ and using theorem 1.2 we may assume that
$x_1,\dots ,x_m$ is a system of parameters of $M/\goth pM$. So we have $\dim
M/\goth pM$=$\dim (R/\goth p)=m$. Therefore we get
$$\eqalign {\Hom _R(R/\goth p,D(\LCMo ^m_I(M)))&=D(\LCMo ^m_I(M/\goth pM))\cr
&=D(\LCMo ^m_\goth m(M/\goth pM)))\cr &\neq 0.\cr }$$
On the other hand we have $\Hom _R(R/\goth q,D(\LCMo ^m_I(M)))=0$ for every
prime ideal $\goth q$ of $R$ containing $\goth p$ properly by (i); the statement
follows.
\par
(iii) Using (ii) it remains to show $x\not\in \goth
p$ for every $\goth p\in \Ass _R(D(\LCMo ^1_{xR}(R)))$. As we have seen our hypothesis
implies $\LCMo ^1_{xR}(R/\goth p)\neq 0$. So we must have $x\not\in \goth p$.
\par
(iv) The first statement follows from (ii) and then the second statement from
(i).
\par
(v) Let $\goth p\in \Supp _R(M)$ such that $\dim(R/\goth p)=m$ and
$\goth p\in \Ass _R(D(\LCMo ^m_I(M)))$. We have to show that
$x_1,\dots ,x_m$ is a system of parameters of $M/\goth pM$: $\LCMo
^m_I(M/\goth pM)\neq 0$ implies $\LCMo ^m_I(R/\goth p)\neq 0$. As
$R$ and hence $R/\goth p$ are complete we may conclude from
Hartshorne-Lichtenbaum vanishing that $\dim (R/(I+\goth p))=0$, i.
e. $x_1,\dots ,x_m$ is a system of parameters of $R/\goth p$.
\bigskip
{\bf 1.4 Remarks}
\par
(i) In theorem 1.3 (ii) we do not have equality in general. For a
counterexample see [11, 2.2.4].
\par
(ii) Statements (ii) and (iii) of theorem 1.3 are generalizations of
[11, theorems 2.2.1 and 2.2.3] resp. of [11, corollary 2.2.2], where
the same statements were shown in special cases depending on the
characteristics of $R$ and $R/\goth m$.
\par
(iii) Obviously the statements of theorem 1.3 fit together well with conjectures (*)
and (+), but are not sufficient to prove either of them.
\bigskip
{\bf 1.5 Lemma}
\par
Let $(S,\goth m)$ be a noetherian local complete Gorenstein ring of
dimension $n+1$ ($\geq 1$) and $\goth P\subseteq S$ a prime ideal of height $n$. Then
$$D(\LCMo ^n_\goth P(S))=\widehat {S_\goth P}/S$$
holds canonically.
\par
Proof: This is a special case of [6, Lemma 3.1].
\bigskip
{\bf 1.6 Lemma}
\par
Let $(R,\goth m)$ be a noetherian local complete domain and
$I\subseteq R$ a prime ideal such that $\dim (R/I)=1$. Then there
exist a noetherian local complete regular ring $S$, a local
homomorphism $S\buildrel \rho \over \to R$ and a prime ideal $\goth
Q\subseteq R$ such that $R$ is finite as an $S$-module and such that
$$\height (\ker (\rho ))=1, \dim (S/\goth Q)=1, \sqrt {\goth QR}=I, \ker
(\rho )\subseteq \goth Q$$ hold. Furthermore, $R$ and $S$ have the
same residue field via $\rho $.
\par
Proof: See [3, Proof of (1.2)]. \bigskip \mittel 2. Results \normal
\bigskip
In this section we calculate the set of associated primes of $D(\LCMo
^{d-1}_J(R))$, where $J$ is a one-dimensional ideal in a $d$-dimensional local
complete ring $R$; furthermore we will partially calculate this set in the more
general situation where $R$ is not necessarily complete.
\bigskip
{\bf 2.1 Lemma}
\par
Let $R$ be a noetherian ring.
\par
(i) Let $\goth P$ be a prime ideal of $R$ which is not maximal. Then the
equivalence
$$R_\goth P=\widehat {R_\goth P}\iff \goth P\hbox { is minimal in }\Spec(R)$$
holds.
\par
(ii) Assume that $R$ is local (and noetherian) and that all prime ideals associated
to $R$ are minimal in $\Spec (R)$. Then $\Ass _R(\hat R/R)\subseteq \Ass (R)$
holds.
In particular if $R$ is a non-complete (local) domain (i. e. if $R\subsetneq
\hat R$),
$$\Ass _R(\hat R/R)=\{ 0\} $$
holds.
\par
Proof: (i) The implication $\Leftarrow $ is clear as every zero-dimensional
local noetherian ring is complete. We assume there exists a prime ideal $P$
of $R$ which is neither minimal nor maximal in $\Spec(R)$ and such that
$R_P=\widehat {R_P}$. $PR_P$ is not minimal in $\Spec
(R)$. Choose $Q,Q^\prime \in \Spec (R)$ such that $Q^\prime
\subsetneq P\subsetneq Q$ and such that $\dim (R_Q/PR_Q)=1$. We set ${\frak
R}:=R_Q/Q^\prime R_Q$ and
${\frak P}:=PR_Q/Q^\prime R_Q\in \Spec ({\frak R})$ and we get
$${\frak R_P}=R_P/Q^\prime R_P=\widehat {R_P}/Q^\prime \widehat {R_P}=\widehat
{R_P/Q^\prime R_P}=\widehat {\frak R_P}.$$
So we may assume that $R$ is a local domain  and $\dim (R/P)=1$.
\par
Take $y\in \goth m\setminus P$. Assume that for some $n\in \Naturalsign $
$$P^{(n)}\subseteq P^{(n+1)}+yR$$
holds ($P^{(n)}:=P^nR_P\cap R$ is a $P$-primary ideal of $R$ such that
$P^{(n)}R_P=P^nR_P$). It would follow that
$$P^{(n)}=P^{(n)}\cap (P^{(n+1)}+yR)=P^{(n+1)}+(P^{(n)}\cap
yR)=P^{(n+1)}+yP^{(n)}$$
and then $P^{(n)}=P^{(n+1)}$ by the Nakayama lemma. Again by Nakayama, this
would imply $P^nR_P=0$ and so $PR_P$ would be minimal in $\Spec (R_P)$. We
conclude that for every $n\in \Naturalsign $
$$P^{(n)}\not\subseteq P^{(n+1)}+yR$$
holds. For every $n\in \Naturalsign $ we choose $x_n\in P^{(n)}\setminus
(P^{(n+1)}+yR)$ and define (for every $n\in \Naturalsign ^+$)
$$\xi _n:=\sum _{i=0}^{n-1}{x_i\over y^{(i+1)^2}}\in R_P.$$
Because of $\xi  _{n+1}-\xi _n={x_n\over y^{(n+1)^2}}\in P^nR_P$ (for every
$n$), we have
$$(\xi _n+P^nR_P)_{n\in \Naturalsign ^+}\in \widehat {R_P}=R_P.$$
There exists a $\xi \in R_P$ such that
$$(\xi +P^nR_P)_{n\in \Naturalsign ^+}=(\xi _n+P^nR_P)_{n\in \Naturalsign
^+},$$
i. e.
$$\xi -\xi _n\in P^nR_P$$ for all $n\in \Naturalsign ^+$.
\par
Write $\xi ={a\over b}$, where $a\in R,b\in R\setminus P$. The ideal $P+bR$ of
$R$ is either $R$ or $\goth m$-primary, so there exist $p\in \Naturalsign ^+$
and $c\in R$ such that $y^p-bc\in P$; it follows that
$$y^{pn}-bc_n\in P^n,$$
where
$$c_n:=b^{-1}(y^{pn}-(y^p-bc)^n)\in R$$
and we conclude that
$$\xi -{ac_n\over y^{pn}}={ay^{pn}-abc_n\over by^{pn}}={a(y^p-bc)^n\over
by^{pn}}\in P^nR_P$$
for every $n\in \Naturalsign ^+$. We get
$$\xi _n-{ac_n\over y^{pn}}=\xi -{ac_n\over y^{pn}}-(\xi -\xi _n)\in P^nR_P$$
for every $n\in \Naturalsign ^+$. From this we get (for $n>p$) after
multiplication by $y^{n^2}$ that
$$\sum _{i=0}^{n-1}x_iy^{n^2-(i+1)^2}-ac_ny^{n(n-p)}\in P^{(n)}$$
and in particular $x_{n-1}\in P^{(n)}+yR$ which is a contradiction.
\par
(ii) We have to prove only the first statement, the second one follows from it
immediately; Let $P$ be an arbitrary element of $\Spec (R)\setminus \Ass (R)$;
We conclude $\Hom _R(R/P,R)=0$ and hence also $\Hom _R(R/P,\hat R)=0$
(because $P$ contains an element which operates injectively on $R$ and $\hat
R$ is flat over $R$). Thus the short exact sequence
$$0\to R\buildrel \subseteq \over \to \hat R\to \hat R/R\to 0$$
induces an exact sequence
$$0\to \Hom _R(R/P,\hat R/R)\to \Ext ^1_R(R/P,R)\buildrel \varphi \over \to
\Ext ^1_R(R/P,\hat R).$$
By our hypothesis there exists $x\in P$ such that $x\not\in Q$ for all $Q\in
\Ass (R)$. We get short exact sequences
$$0\to R\buildrel x\over \to R\to R/xR\to 0$$
and
$$0\to \hat R\buildrel x\over \to \hat R\to \hat R/x\hat R\to 0.$$
Because of $x\in P$ a commutative diagram with exact rows is induced:
$$\matrix{&0&\to &\Hom _R(R/P,R/xR)&\to &\Ext ^1_R(R/P,R)&\to &0\cr
&&&\downarrow \psi &&\downarrow \varphi &&\cr &0&\to &\Hom _R(R/P,\hat R/x\hat
R)&\to &\Ext ^1_R(R/P,\hat R)&\to &0.\cr }$$
$\psi $ is injective as $R/xR\subseteq \widehat {R/xR}=\hat R/x\hat R$.
Therefore $\varphi $ is injective which implies that $\Hom _R(R/P,\hat
R/R)=0$, i. e. $P\not\in \Ass _R(\hat R/R)$.
\bigskip
{\bf 2.2 Theorem}
\par
Let $(R,\goth m)$ be a $d$-dimensional local noetherian complete
domain, where $d\geq 1$; let $P$ be a prime ideal of $R$ such that
$\dim (R/P)=1$. Then
$$\{ 0\} \in \Ass _R(D(\LCMo ^{d-1}_P(R)))$$
holds.
\par
Proof: We apply lemma 1.6, set $R_0:=\im (\rho )$ and consider the
ideal $\goth Q$ from 1.6 as an ideal of $R_0$. By lemma 1.6, $R_0$
is a complete intersection, in particular it is Gorenstein. By
$\goth m_0$ we denote the maximal ideal of $R_0$. $R$ is finite over
$R_0$ and so we have $D_{R_0}(R)=\Hom _{R_0}(R,\InjH
_{R_0}(R_0/\goth m_0)))=\InjH _R(R/\goth m)=D_R(R)$, which implies
$D_{R_0}(M)=D_R(M)$ for every $R$-module $M$. On the other hand the
functor $\LCMo ^{d-1}_\goth Q(\_ )$ is right exact by
Hartshorne-Lichtenbaum vanishing; in particular we have
$$D_R(\LCMo ^{d-1}_\goth P(R))=D_{R_0}(\LCMo ^{d-1}_\goth Q(R_0)\otimes
_{R_0}R)=\Hom _{R_0}(R,D_{R_0}(\LCMo ^{d-1}_\goth Q(R_0)))$$ and so
every $R_0$-monomorphism $R_0\to D_{R_0}(\LCMo ^{d-1}_\goth Q(R_0))$
induces an $R$-monomorphism $\Hom _{R_0}(R,R_0)\to D_R(\LCMo
^{d-1}_\goth P(R))$. But $\{ 0\} \in \Ass _{R_0}(\Hom
_{R_0}(R,R_0))$ holds (because
$$\Hom _{R_0}(R,R_0)\otimes
_{R_0}Q(R_0)=\Hom _{Q(R_0)}(R\otimes _{R_0}Q(R_0),Q(R_0))$$ is a
non-zero $Q(R_0)-$vector space, here we use the fact that $R$ is
finite over $R_0$) and thus it suffices to show $\{ 0\} \in \Ass
_{R_0}(D_{R_0}(\LCMo ^{d-1}_\goth Q(R_0)))$, i. e. we may assume
$R_0$ is Gorenstein. Now, by lemma 1.5, we have a commutative
diagram with exact rows:
$$\matrix {&0&\to &R&\buildrel \subseteq \over \to &\widehat {R_P}&\to &D(\LCMo
^{d-1}_P(R))&\to &0\cr &&&\downarrow \subseteq &&\downarrow = &&&&\cr &0&\to
&R_P&\buildrel \subseteq \over \to &\widehat {R_P}&\to &\widehat {R_P}/R_p&\to
&0.\cr }$$
This diagram induces an epimorphism
$$D(\LCMo ^{d-1}_P(R))\to \widehat {R_P}/R_P.$$
By lemma 2.1 (i) we have $\widehat {R_P}/R_P\neq 0$ and it follows
from lemma 2.1 (ii) that $(\widehat {R_P}/R_P)\otimes _R(Q(R))\neq
0$. Thus we have $D(\LCMo ^{d-1}_P(R))\otimes _RQ(R)\neq 0$ by the
above epimorphism.
\bigskip
{\bf 2.3 Lemma}
\par
Let $(R,\goth m)$ be a noetherian local complete ring, $d:=\dim
(R)\geq 1$, $J\subseteq R$ an ideal of $R$ such that $\dim (R/J)=1$.
Then
$$\{ Q\in \Ass _R(D(\LCMo ^{d-1}_J(R)))\vert \dim (R/Q)=d\} =\{ Q\in \Assh
(R)\vert \dim (R/(J+Q))\geq 1\} .$$ If $\LCMo ^d_J(R)=0$,
$$\Assh (D(\LCMo ^{d-1}_J(R)))=\Assh (R).$$
Proof: The second statement follows from the first one by
Hartshorne-Lichtenbaum vanishing. We prove the first statement now in the
special case where $\dim (R/(J+Q))\geq 1$ for all $Q\in \Assh (R)$. In the
second part of the proof we will show how to reduce the general to the special
situation. Now, in the special case
it suffices to show the inclusion ``$\supseteq $''. By Hartshorne-Lichtenbaum
vanishing we have $\LCMo ^d_J(R)=0$, i. e. $\LCMo ^{d-1}_J(\_ )$ is right
exact. Let $Q\in \Assh (R)$ be arbitrary. The canonical epimorphism $R\to
R/Q=:\overline R$ induces a monomorphism
$$D_{\overline R}(\LCMo ^{d-1}_{J\overline R}(\overline R))=D_R(\LCMo ^{d-1}_J(\overline
R))\to D(\LCMo ^{d-1}_J(R)).$$ If $\{ 0\} \in \Ass _{\overline
R}(D_{\overline R}(\LCMo ^{d-1}_{J\overline R}(\overline R)))$ then
$Q\in \Ass _R(D(\LCMo ^{d-1}_J(R)))$, and so we may assume that $R$
is a domain. If we can write $J=J_1\cap J_2$ with ideals $J_1,J_2$
of $R$ such that $J_1+J_2$ is $\goth m$-primary then because of
$\LCMo ^d_{J_1}(R)=\LCMo ^d_{J_2}(R)=0$ (Hartshorne-Lichtenbaum
vanishing) a Mayer-Vietoris sequence argument gives us an
epimorphism
$$\LCMo ^{d-1}_J(R)\to \LCMo ^d_{J_1+J_2}(R)=\LCMo ^d_\goth m(R).$$
But then theorem 1.3 shows that
$$\{ 0\} =\Assh (R)=\Ass _R(D(\LCMo ^d_\goth m(R)))\subseteq \Ass _R(D(\LCMo
^{d-1}_J(R))).$$ If there is no such decomposition $J=J_1\cap J_2$
of $J$ we may assume that $J$ is a prime ideal; but then the
statement follows from theorem 2.2. Now we turn to the general case,
i. e. we assume there is a $Q\in \Assh (R)$ such that $\dim
(R/(J+Q))=0$. We define $U(R)$ to be the intersection of all
$Q^\prime $-primary components of a primary decomposition of the
zero ideal in $R$ for all $Q^\prime \in \Assh (R)$. Apparantly we
have $\Ass _R(R/U(R))=\Assh (R)$ and $\dim (U(R))<d$. Because of the
latter fact the short exact sequence $0\to U(R)\buildrel \subseteq
\over \to R\to R/U(R)\to 0$ induces an exact sequence
$$0\to D(\LCMo ^{d-1}_J(R/U(R)))\to D(\LCMo ^{d-1}_J(R))\to D(\LCMo
^{d-1}_J(U(R))).$$ Trivially $\dim _R(\Supp _R(\LCMo
^{d-1}_J(U(R))))\leq d-1$ holds. By considering $R/U(R)$ rather then
$R$ we may assume that $\Ass _R(R)=\Assh (R)$. We write $0=I^\prime
\cap I^{\prime \prime }$ with ideals $I^\prime ,I^{\prime \prime }$
of $R$ such that $\Ass _R(R)=\Ass _R(R/I^\prime )\cup \Ass
_R(R/I^{\prime \prime })$ and $\dim (R/(J+Q^\prime ))\geq 1$ for all
$Q^\prime \in \Ass _R(R/I^\prime )$ and $\dim (R/(J+Q^{\prime \prime
}))=0$ for all $Q^{\prime \prime }\in \Ass _R(R/I^{\prime \prime
}))$. It follows that $\dim (R/(J+I^{\prime \prime }))=0$. By using
a Mayer-Vietoris argument and the facts that $\LCMo ^d_J(R/I^\prime
)=0$ (Hartshorne-Lichtenbaum vanishing) and $\LCMo ^i_J(R/I^{\prime
\prime })=\LCMo ^i_\goth m(R/I^{\prime \prime })$ for all $i\in
\Naturalsign $ we get a short exact sequence
$$D(\LCMo ^{d-1}_\goth m(R/I^\prime +I^{\prime \prime }))\to D(\LCMo
^{d-1}_J(R/I^\prime ))\oplus D(\LCMo ^{d-1}_\goth m(R/I^{\prime \prime }))\to
$$
$$\to D(\LCMo ^{d-1}_J(R))\to D(\LCMo ^{d-2}_\goth m(R/(I^\prime +I^{\prime
\prime }))).$$
It is clear that we have
$$\dim _R(\Supp _R(D(\LCMo ^{d-1}_\goth m(R/(I^\prime +I^{\prime \prime
}))))),\dim _R(\Supp _R(D(\LCMo ^{d-2}_\goth m(R/(I^\prime
+I^{\prime \prime })))))\leq d-1.$$ One can write $R$ as a quotient
of a local Gorenstein ring $S$ such that $\dim(S)=\dim(R)$; over $S$
one has local duality and therefore one may conclude (note that in
the formula below it does not make any difference if we take $D$
with respect to $R$ or to $S$)
$$D(\LCMo ^{d-1}_\goth m(R/I^{\prime \prime }))=\Ext
^1_S(R/I^{\prime \prime },S)$$ and hence
$$\dim _R(\Supp _R(D(\LCMo ^{d-1}_\goth m(R/I^{\prime \prime }))))\leq d-1.$$
Thus we get, by what is already shown,
$$\{ Q\in \Ass _R(D(\LCMo ^{d-1}_J(R)))\vert \dim (R/Q)=d\} =\{ Q\in \Ass
_R(D(\LCMo ^{d-1}_J(R/I^\prime )))\vert \dim (R/Q)=d\} =\Assh (R/I^\prime ).$$
\bigskip
{\bf 2.4 Theorem}
\par
Let $(R,\goth m)$ be a $d$-dimensional noetherian local ring
and $J\subseteq R$ an ideal such that $\dim (R/J)=1$ and $\LCMo
^d_J(R)=0$. Then
$$\Assh (R)=\Assh (D(\LCMo ^{d-1}_J(R))).$$
Proof: One has $\LCMo ^d_{J\hat R}(\hat R)=\LCMo ^d_J(R)\otimes _R\hat R=0$
and
$$\eqalign {D_{\hat R}(\LCMo ^{d-1}_{J\hat R}(\hat R))&=D_{\hat R}(\LCMo
^{d-1}_J(R)\otimes \hat R)\cr &=\Hom _R(\LCMo ^{d-1}_J(R),D_{\hat R}(\hat
R))\cr &=D_R(\LCMo ^{d-1}_J(R))\cr }$$
These imply that every $\hat R$-monomorphism $\varphi :\hat R/\goth P\to
D_{\hat R}(\LCMo ^{d-1}_{J\hat R}(\hat R))$, where $\goth P$ is a prime ideal
of $\hat R$, induces an $R$-monomorphism $\hat R/\goth P\to D_R(\LCMo
^{d-1}_J(R))$. On the other hand we have a $R$-monomorphism $R/\goth p\to \hat
R/\goth P$, where $\goth p:=\goth P\cap R$. Composition of these monomorphisms
gives us a monomorphism
$$R/\goth p\to D_R(\LCMo ^{d-1}_J(R)).$$
Because of $\Assh (R)=\{ \goth P\cap R\vert \goth P\in \Assh (\hat
R)\}$ we may assume $R$ is complete. But then the statement follows
from lemma 2.3.
\bigskip
{\bf 2.5 Theorem}
\par
Let $R$ be a $d$-dimensional local complete ring and
$J\subseteq R$ an ideal such that $\dim (R/J)=1$ and $\LCMo ^d_J(R)=0$. Then
$$\Ass _R(D(\LCMo ^{d-1}_J(R))=\{ P\in \Spec(R)\vert \dim (R/P)=d-1,\dim
(R/(P+J))=0\} \cup \Assh (R).$$
Proof: Let $P\in \Spec (R)$. If $\dim (R/P)\leq d-2$ we have $\Hom
_R(R/P,D(\LCMo ^{d-1}_J(R)))=D(\LCMo ^{d-1}_J(R/P))=0$
and hence $P\not\in \Ass _R(D(\LCMo^ {d-1}_J(R)))$. If $\dim (R/P)=1$ then
($\overline R:=R/P$):
$$\Hom _R(R/P,D(\LCMo ^{d-1}_J(R)))=D(\LCMo ^{d-1}_J(R/P))=D(\LCMo
^{d-1}_{J\overline R}(\overline R)).$$ $R$ is complete and so, by
Hartshorne-Lichtenbaum the equivalence
$$\LCMo ^{d-1}_{J\overline R}(\overline R)\neq 0\iff \dim (\overline R/J\overline R)=0$$
holds. On the other hand we have $\overline R/J\overline R=R/(P+J)$ and
therefore we have
$$\{ P\in \Ass _R(D(\LCMo ^{d-1}_J(R)))\vert \dim (R/P)=d-1\}=$$
$$=\{ P\in \Spec (R)\vert \dim (R/P)=d-1,\dim (R/(P+J))=0.$$
Now the statement follows from lemma 2.3.
\bigskip
{\bf References}
\normal
\smallskip
\parindent=0.8cm
\def\litem{\par\noindent \hangindent=\parindent\ltextindent}
\def\ltextindent#1{\hbox to \hangindent{#1\hss}\ignorespaces}
\litem{1.} Bass, H. On the ubiquity of Gorenstein rings, {\it Math. Z.} {\bf
82}, (1963) 8--28.
\medskip
\litem{2.} Brodmann, M. and Hellus, M. Cohomological patterns of coherent
sheaves over projective schemes, {\it Journal of Pure and Applied Algebra}
{\bf 172}, (2002) 165--182.
\medskip
\litem{3.} Brodmann, M. and Huneke, C. A quick proof of the
Hartshorne-Lichtenbaum vanishing theorem, Algebraic geometry and its
applications, {\it Springer, New York} , (1994) 305--308.
\medskip
\litem{4.} Brodmann, M. and Sharp, R. J. Local Cohomology, {\it
Cambridge studies in advanced mathematics} {\bf 60}, (1998).
\medskip
\litem{5.} Bruns, W. and Herzog, J. Cohen-Macaulay Rings, {\it
Cambridge University Press}, (1993).
\medskip
\litem{6.} Chiriacescu, G. Cofiniteness of local cohomology modules
over regular local rings, {\it Math. Soc.} {\bf 32}, (2000) 1--7.
\medskip
\litem{7.} Eisenbud, D. Commutative Algebra with A View Toward
Algebraic Geometry, {\it Springer Verlag}, (1995).
\medskip
\litem{8.} Grothendieck, A. Local Cohomology, {\it Lecture Notes in
Mathematics, Springer Verlag}, (1967).
\medskip
\litem{9.} Hellus, M. Matlis duals of top local cohomology modules
and the arithmetic rank of an ideal, to appear in {\it
Communications in Algebra}.
\medskip
\litem{10.} Hellus, M. On the set of associated primes of a local
cohomology module, {\it J. Algebra} {\bf 237}, (2001) 406--419.
\medskip
\litem{11.} Hellus, M. On the associated primes of Matlis duals of
top local cohomology modules, to appear in {\it Communications in
Algebra} {\bf 33}.
\medskip
\litem{12.} Huneke, C. Problems on Local Cohomology, {\it Res. Notes
Math. } {\bf 2}, (1992) 93--108.
\medskip
\litem{13.} Huneke, C. and Lyubeznik, G. On the vanishing of local
cohomology modules, {\it Invent. math.} {\bf 102}, (1990) 73-93.
\medskip
\litem{14.} Matlis, E. Injective modules over Noetherian rings, {\it
Pacific J. Math.} {\bf 8}, (1958) 511--528.
\medskip
\litem{15.} Matsumura, H. Commutative ring theory, {\it Cambridge
University Press}, (1986).
\medskip
\litem{16.} Scheja, G. and Storch, U. Regular Sequences and
Resultants, {\it AK Peters}, (2001).
\end